\documentclass{article}
\usepackage{amsmath}
\usepackage{amssymb}
\usepackage{hhline}

\newtheorem{theorem}{Theorem}
\newtheorem{theoremc}{Theorem}
\newtheorem{theoremd}{Theorem}
\newtheorem{theoreme}{Theorem}
\newtheorem{theoremr}{Theorem}

\newtheorem{cor}[theoremc]{Corollary}
\newtheorem{dfn}[theoremd]{Definition.}
\newtheorem{examp}[theoreme]{Example}
\newtheorem{lem}[theorem]{Lemma}
\newtheorem{prop}[theorem]{Proposition}
\newtheorem{rk}[theoremr]{Remark}
\newenvironment{proof}[1][Proof]{\textbf{#1.} }{\qed}

\newcommand\bib[1]{\bibitem[#1]{#1}}

\newcommand\1{{\bf 1}}
\renewcommand\a{\alpha}
\renewcommand\b{\beta}
\newcommand\bff[1]{\vspace{4pt}\underline{\bf #1}$^\circ$)}

\newcommand\C{{\mathcal C}}
\newcommand\CC{{\mathbb C}}
\newcommand\D{{\cal D}}
\newcommand\e{\eta}

\newcommand\g{\gamma}
\newcommand\G{{\cal G}}
\renewcommand\H{{\cal H}}
\newcommand\hht{h_{\text{top}}}
\newcommand\La{\Lambda}
\renewcommand\l{\lambda}
\newcommand\m{\mu}

\renewcommand\O{\Omega}
\newcommand\oo{\omega}
\newcommand\op[1]{\mathop{\rm #1}\nolimits}

\newcommand\p{\partial}

\newcommand\po{\hspace{-4pt}{\bf .}\,\,}

\newcommand\qed{\phantom{\underline{y}}\hfill\hfill$\square$}

\newcommand\R{{\mathbb R}}
\renewcommand\t{\times}

\newcommand\Te{\Theta}
\newcommand\te{\theta}

\newcommand\ve{\varepsilon}
\newcommand\vp{\varphi}
\newcommand\we{\wedge}
\newcommand\x{\xi}
\newcommand\z{\sigma}
\newcommand\Z{{\mathbb Z}}

\begin{document}
 \title{Examples of integrable \\ sub-Riemannian geodesic flows}
 \author{Boris Kruglikov%
 }
 \date{}
\maketitle

 \begin{abstract}
We exhibit examples of sub-Riemannian metrics with integrable geodesic
flows and positive topological entropy.
 \end{abstract}

\section*{Introduction}

\hspace{13.5pt}
Consider a distribution on a manifold $M^m$, i.e.\ subbundle of the tangent
bundle $\Pi\subset TM$. Non-holonomic Riemannian metric is a Riemannian
metric $g\in S^2\Pi^*$ on this bundle. We call the pair $(\Pi,g)$
{\em sub-Riemannian structure\/}. A curve $\g:[0,1]\to M$ is called
horizontal if $\dot\g$ is a section of $\Pi$. We denote the space of
horizontal curves joining $x$ to $y$ by $\H(x,y)$. A theorem of
Rashevsky-Chow (\cite{R}) states that if $M$ is connected and $\Pi$ is
completely non-holonomic then $\H(x,y)$ is always non-empty. By completely
non-holonomic we mean distribution $\Pi$, such that the module
$\D_\Pi^{(N)}$ of order $\le N$ self-commutators (of various kinds)
of sections of $\Pi$ is equal to the
module $\D(M)$ of all vector fields for some big $N$.
From now on we consider only completely non-holonomic distributions.

For horizontal curves we calculate its length
$l_g(\g)=\int_0^1\|\dot\g\|_gdt$ and this produces
{\em sub-Riemannian\/} distance (metric) on $M$ by
 $$
d_g(x,y)=\inf\limits_{\g\in\H(x,y)}l_g(\g).
 $$
A curve $\g\in\H$ is called {\em geodesic} if it realizes the minimum
sub-distance for any two of its close points. The description of the most
geodesics (normal ones) is given by the Euler-Lagrange variational
principle. There is a Hamiltonian reformulation of this principle, due to
Pontrjagin and co-authors \cite{PBGM}, which allows to consider
the geodesic flow as the usual Hamiltonian flow on $T^*M$. There appear
occasionally geodesics of different kind -- abnormals -- which are not
governed by the Pontrjagin principle for $\g$, but depend on the
distribution $\Pi$ only. However if we consider contact distributions
$\Pi$, i.e.\ distributions such that for any non-zero section $\a$ of the
bundle $\op{Ann}(\Pi)\subset T^*M$ we have $\a\we(d\a)^n\ne0$ for $m=2n+1$
(in particular $m=\op{dim}M$ is odd), then all geodesics
are normal.

As in the standard theory of geodesics we say the metric
$g$ is integrable if the Hamiltonian flow of this metric
is integrable on $T^*M$ in the Liouville sense, i.e.\ there are
a.e.\ functionally independent integrals $I_1=H,I_2,\dots,I_n$ which
Poisson-commute $\{I_k,I_l\}=0$ (see \cite{BF} for a review of
methods and problems).

For any diffeomorphism $f:M\to M$ the suspension $M\tilde\t_fS^1$ is the
quotient of $M\t\R$ by the equivalence relation $(x,t)\sim(f(x),t+1)$. In
particular for any matrix $A\in\op{SL}_n(\Z)$ we have a linear automorphism
$A:T^n\to T^n$ of the torus and hence the manifold $M_A=T^n\tilde\t_AS^1$.
In the paper \cite{BT1} a series of examples of (Riemannian) geodesic flows
on $M_A$ was constructed with the property that they are
$C^\infty$-integrable, $C^\oo$- (analytically) nonintegrable
and have the topological entropy $h_{top}>0$.

We show that similar effect takes place in sub-Riemannian geometry too.

 \begin{theorem}\po\label{th1}
For every $A\in\op{SL}_2(\Z)$ there are a contact distribution $\Pi$ on
$M_A=T^2\tilde\t_AS^1$ and a non-holonomic Riemannian structure $g$ on it
such that the sub-Riemannian geodesic flow is $C^\infty$-integrable.

In the case of real eigenvalues for $A\ne\pm E$ the flow is
not integrable with geometrically simple set of integrals. Moreover
the topological entropy is positive in semi-simple non-trivial case,
$\hht=\max\limits_{\l_i\in\op{Sp}(A)}\{\ln|\l_1|,\ln|\l_2|\}$,
and vanish otherwise, $\hht=0$.
 \end{theorem}

The geometrically simple set of integrals includes analytic polynomial by
momenta integrals and is defined in \S\ref{S4} (remark that we consider
only integrability with Liouville tori, not cylinders, see
\S\ref{S1}\!\! \ref{S1.1}).

Note that integrability of sub-Riemannian flow is a more complicated
(and fascinating) fact since the description of the geodesics is
completely nontrivial even locally (see \cite{VG} for a picture
of sub-Riemannian wave front).

In particular locally over a point $x$ of the manifold
Riemannian geodesic flow in $\tau_M^{-1}(U_x)$ is integrable,
while in non-holonomic case it is no longer true. On the other hand
integrability does not imply good behaviour in other
respects: The well-known Martinet case (see ex.\ref{Martinet}) is
$C^\oo$-integrable, but small geodesic balls are not sub-analytic
(\cite{AS}) and there are abnormal geodesics.

 \begin{rk}\po
One can easily generalize the above examples to the examples of integrable
sub-Riemannian geodesic flows on $M_A$ for arbitrary $A\in\op{SL}_n(\Z)$
with $n>2$. There appear distributions of various ranks and kinds (they
correspond to higher Vandermonde determinants instead of $2\t2$ matrices
given below).  Though completely non-holonomic, they have very degenerate
curvature, i.e.\ the map $\Te_\Pi:\La^2\Pi\to\nu=TM/\Pi$ induced by
commutators (in particular, the growth vector can be rather long). They
never produce higher-dimensional contact examples.
 \end{rk}

Every sub-Riemannian metric $g$ on a manifold of
$\dim=2n+1$ produces canonically a transversal to the contact distribution
$\Pi$ vector field -- Reeb field $\nu_\a$ (also denoted $R_\a$). To see it
we choose $\a\in\op{Ann}\Pi$ such that $\|\left.d\a\right|_\Pi\|=1$, where
the norm is induced by the metric. In other words
$\left.(d\a)^n\right|_\Pi=vol_g$.  This gives $\a$ uniquely for odd $n$
and up to $\pm1$ for even if the distribution $\Pi$ is oriented. The Reeb
field is defined now uniquely by $\a(\nu_\a)=1$, $i_{\nu_\a}d\a=0$.

 \begin{theorem}\po\label{th2}
In addition within the construction of theorem~\ref{th1} if\/ $\Pi$ is
cooriented {\rm(}$\op{Sp}A$ is contained in either $\R_+$ or
$S^1\subset\CC${\rm)} the topological entropy of the Reeb flow is
$h_{top}(\nu_\a)=0$.
 \end{theorem}

We introduce Riemannian metric $\hat g$ on $M$ by the requirements
$\left.\hat g\right|_\Pi=g$, $\nu_\a\perp_{\hat g}\Pi$,
$\|\nu_\a\|_{\hat g}=1$.

 \begin{theorem}\po\label{th3}
The geodesic flow of Riemannian metric $\hat g$ on $M_A$ is
$C^\infty$-integrable, $C^\oo$-nonintegrable for real eigenvalues of
$A\ne\pm E$ and has positive topological entropy for semi-simple
nontrivial case.
 \end{theorem}

In fact the metric $\hat g$ coincides with the metric on $M_A$ introduced
by Bolsinov and Taimanov. The first example of $C^\infty$-integrable
geodesic flow violating obstructions for
$C^\oo$-integrability was found by Butler \cite{B} and then it was
included into a much bigger family of examples-suspensions in~\cite{BT1}.
Their examples obviously generalize to Lorentian and other
semi-geometries. We show the case of sub-Riemannian geometry (actually of
any dim) can be considered as well.

 \begin{rk}\po
Note that examples from {\rm\cite{BT1}} are easily
generalized to a larger group
$\op{GL}_n(\Z)=\op{SL}_n(\Z)\t\Z_2$. However the sub-Riemannian case
is different: the manifold $M_A$ with $\det A=-1$ does not admit any
contact structure. Actually a contact structure on a 3-dimensional
manifold gives a canonical orientation.
 \end{rk}

We discuss some relations between $g$, $\hat g$ and $\nu_\a$
at the end of the paper. In \S~\ref{S6} we exibit Poisson action of
$\R^3$ such that entropy of every non-zero vector is non-zero.

Let us note that theorem~\ref{th1} provides the first example of
$C^\oo$-polynomial in momenta nonintegrable sub-Riemannian geodesic flow.
In paper~\cite{MSS} it was presented only algebraically non-integrable
sub-Riemannian flow on a Lie group of $\dim M=6$. In dimension 3
all left-invariant sub-geodesic flows on Lie groups
are Liouville integrable (as we show in \S\ref{S3}).

{\footnotesize Acknowledgements.
I would like to thank J.\,P.\ Gauthier, A.\,V.\ Bolsinov, V.\,S.\ Matveev
and O.\ Kozlovski for helpful discussions. I am indebted to A.\,B.\ Katok
for the reference in \S\ref{commut}. My interest to sub-Riemannian geometry
has grown up during a Singularity Theory workshop at Newton Institute
(Cambridge), 2000. I would like to thank the Institute and especially
V.\,V.\ Goryunov and V.\,M.\ Zakalyukin for the organization of this
activity and hospitality.}

\section{Sub-Riemannian geometry and Hamiltonian systems}\label{S1}
\subsection{\sl Integrable Hamiltonian flows.}\label{S1.1}

\hspace{13.5pt}
Let $(W^{2n},\O)$ be a symplectic manifold. Hamiltonian vector field is
the field $\O$-dual to an exact 1-form $dH$: $\O(X_H,\cdot)=dH$.
We denote the field $X_H$ also by $\op{sgrad}H$, since the above is
similar to the definition of the usual gradient (and in order not to
confuse with $M_A$ above).

Poisson brackets can be defined as $\{F,G\}=\O(X_F,X_G)$. The Hamiltonian
system is called complete (or Liouville) integrable if there are
additionally to $I_1=H$ involutive integrals $I_2,\dots,I_n$,
$\{I_j,I_k\}=0$, which are functionally independent a.e. By Liouville
theorem (\cite{A}) a full measure set $W'\subset W$ is then foliated by
cylinders (tori in compact case), and each cylinder has a neighborhood
with coordinates $\vp\in T^{n-r}\t\R^r$ and new $I\in\R^n$ s.t.\
$\O=dI\we d\vp$, $\{I_j=c_j\}_{1\le j\le n}\simeq T^{n-r}\t\R^r$ and the
flow is $\dot\vp=\oo(I)$. The number $r$ (noncompactness rank) vanishes in
many important examples. We will always suppose $r=0$ (even for noncompact
$W$) and call this case {\em Liouville integrability\/}.

Let $(M,g)$ be a Riemannian manifold. Then the geodesic flow
$\vp_t:TM\to TM$ can be considered as Hamiltonian on
$T^*M\stackrel{\sharp^g}\simeq TM$ with the standard
symplectic structure $\O$ if we choose the Hamiltonian $H=\frac12\|p\|^2$,
$p\in T^*M$. We say $g$ is integrable if the flow $\vp_t$ is
Liouville-integrable.

\subsection{\sl Pontrjagin maximum principle.}

\hspace{13.5pt}
Consider now non-holonomic case. We start with an arbitrary completely
non-holonomic distribution $\Pi\subset TM$.

Let $(T^*M,\O)$ be the cotangent bundle equipped with the standard
symplectic structure. The non-holonomic metric defines isomorphism
$\sharp^g:\Pi^*\to\Pi$.  Consider the inclusion $i:\Pi\hookrightarrow TM$.
Then we have vector bundles morphism $\Psi_g$ defined as the following
composition:
 $$
T^*M\stackrel{i^*}\to\Pi^*\stackrel{\sharp^g}\to\Pi\stackrel{i}\to TM.
 $$
Contrary to Riemannian situation this is not isomorphism. We have:
$\op{Ker}(\Psi_g)=\op{Ann}\Pi$ and $\op{CoKer}(\Psi_g)=\nu=TM/\Pi$.

Define the Hamiltonian function $H$ on $T^*M$ as the composition
 $$
T^*M\stackrel{i^*}\to\Pi^*\stackrel{\sharp_g}\to\Pi
\stackrel{\frac12\|\cdot\|^2_g}\to\R.
 $$
This function can be locally described as follows. Let $\x_1,\dots,\x_k$
be some orthonormal basis of vector fields tangent to $\Pi$. Every vector
field is a fiber-linear function on $T^*M$. So we have
$H=\frac12\sum_1^k\x_i^2$.

The Pontrjagin maximum principle (\cite{PBGM}) states that trajectories of
this vector fields in region $\{H>0\}$ projected to $M$ are optimal for
the corresponding variational problem. They are called (normal) geodesics.

 \begin{examp}\po\label{ex:xe}
Consider $M=T^3=\R^3/2\pi\Z^3$ with cyclic coordinates
$(\vp_1,\vp_2,\vp_3)$.  Let $\a=\sin\vp_1d\vp_2+\cos\vp_1d\vp_3$ be a
contact form and $\Pi=\op{Ker}\a$ the corresponding distribution. Let the
metric be induced from the standard metric $ds^2=\sum d\vp_i^2$ on the
torus. Then $H=\frac12[p_1^2+(\cos\vp_1p_2-\sin\vp_1p_3)^2]$ and the
Hamilton equations $\dot\vp=\p H/\p p$, $\dot p=-\p H/\p\vp$ have the
form:
 $$
 \begin{aligned}
\dot\vp_1&=p_1,\quad
\dot\vp_2=\cos^2(\vp_1)p_2-\frac{\sin(2\vp_1)}2 p_3,\quad
\dot\vp_3=-\frac{\sin(2\vp_1)}2 p_2+\sin^2(\vp_1)p_3,&\\
\dot p_1&=\frac{\sin(2\vp_1)}2 (p_2^2-p_3^2)+\cos(2\vp_1)p_2p_3,\quad
\dot p_2=0,\quad
\dot p_3=0.&
 \end{aligned}
 $$
They can be easily integrated. The functions $I_2=p_2$, $I_3=p_3$ are
obviously integrals. This can be also checked via Poisson brackets: the
functions $I_1=H$, $I_2$, $I_3$ are involutive: $\{I_j,I_k\}=0$. The only
singularities are $\bigodot$ and $\bigcirc\hspace{-5pt}\bigcirc$
(so-called atoms $A$ and $C_2$).
 \end{examp}

 \begin{examp}\po\label{HSR}
The geodesics on the Heisenberg group are given by the
$2H=(p_1+x_3p_2)^2+p_3^2$. Solving the Hamiltonian equation we see that
geodesics are spirals in the direction of $x_2$-axis, projecting to
arbitrary (including radius $\infty$) circles on the plane $\R^2(x_1,x_3)$.
 \end{examp}

\subsection{\sl Note on abnormals.}

\hspace{13.5pt}
In some situations there are optimal lines, which depend on the
distribution $\Pi$ only. They are called abnormal geodesics and
in generic case are described as follows.

Consider the submanifold $S=\op{Ann}\Pi\subset T^*M$. It is given by
equation $H=0$. So contrary to the Riemannian case
$k=\op{codim}S=\op{dim}\Pi<m=\dim M$.
Let $0_M\subset T^*M$ be the zero section
and $S_0=S\setminus 0_M$ be the complement.
This $S_0$ is a punctured cone over $0_M$.

Let $k\in2\Z$. Then generically there is a hypersurface
$\Sigma^{2n-k-1}\subset S_0$ (maybe empty),
where $\left.\O^{n-k/2}\right|_{S_0}=0$.
The integral lines of $\op{Ker}\left(\left.\O\right|_\Sigma\right)$
projects to the abnormals.

Let $k\in\{2\Z+1\}$. Then generically kernels of $\left.\O\right|_{S_0}$
are one-dimensional, so they integrates to lines projecting to the
abnormals.

This geometric description comes from the variational approach
of~\cite{Hsu}.

 \begin{examp}\label{Martinet}\po
Consider the Martinet case: $\Pi=\op{Ker}\a$ for $\a=dy-z^2dx$ on
$\R^3(x,y,z)$. Let the metric be lifted from the Euclidean plane
$\R^2(x,z)$ (recall the abnormals do not depend on a choice of the
metric). So $2H=(p_x+z^2p_y)^2+p_z^2$. There is no flow on $S=\{H=0\}$
(as usual because $H$ is quadratic). However there are geodesics.

Actually, $S=\{p_x=-z^2p_y,p_z=0\}$. So the coordinates are $(x,y,z,p_y)$.
Since $\O=dx\we dp_x+dy\we dp_y+dz\we dp_z$ we get
$\left.\O^2\right|_S=2zp_ydx\we dy\we dz\we dp_y$. Since on $S_0$ $p_y\ne0$
we conclude that $\Sigma^3=\{z=0\}\subset S_0$ with coordinates
$(x,y,p_y)$. Now $\left.\O\right|_\Sigma=dy\we dp_y$, whence
$\op{Ker}\O=\langle\p_x\rangle$ and the abnormal geodesics are given by
$\{y=\op{const},z=0\}$.
 \end{examp}

 \begin{examp}\po
Consider the Engel distribution:
$\op{Ann}\Pi=\langle dy-zdx,dz-wdx\rangle$, i.e.
$\Pi=\langle \p_x+z\p_y+w\p_z,\p_w\rangle$ on $\R^4(x,y,z,w)$.
The manifold $S$ is given by $\{p_x=-zp_y-wp_z,p_w=0\}$ and has
coordinates $(x,y,z,w,p_y,p_z)$. So $\left.\O^3\right|_{S_0}=0$ is
equivalent to $p_z=0$. This gives the submanifold $\Sigma^5$ and
$\left.\O\right|_\Sigma=-p_ydx\we dz-zdx\we dp_y+dy\we dp_y$, whence
$\op{Ker}(\left.\O\right|_\Sigma)=\langle\p_w\rangle$. This is so-called
characteristic direction for the Engel distribution and the abnormals are
just the corresponding integral curves.
 \end{examp}

It follows from the above description that non-holonomic metrics on
contact distributions have no abnormals.

\subsection{\sl Local integrability.}

\hspace{13.5pt}
Though in many respects sub-Riemannian structures are similar to
Riemannian ones (e.g.\,\cite{ACG}), as the previous subsection shows there
are differences. Another discriminating aspect is Liouville integrability.

Let us consider Euclidean metric locally in a neighborhood $B_\ve(x)$.
Symplectic reduction of the corresponding Hamiltonian flow on $\{H=1/2\}$
is $T^*_{<\ve}S^{n-1}$. By convexity theorem for arbitrary Riemannian
metric and small $\ve$ the symplectic reduction preserves differentiable
type of $T^*S^{n-1}$. Moreover by Moser isotopy method applied to
(open) neighborhood of the zero section we get that germs of reductions
are symplectomorphic in Riemannian and linearized (Euclidean) cases.
So one can find involute set of functions and integrate (by Liouville)
the Riemannian geodesic flow on a neighborhood $U_x\subset B_\ve(x)$.

In contrast even the differentiable type of symplectic reduction for
sub-Riemannian structure in the simplest Heisenberg case can be changed
under arbitrary small perturbation. That is because the geodesics
are curved and being perturbed may leave/return the domain.

\section{Entropy of a dynamical system}\label{S2}
\subsection{\sl Definition of entropies.}

\hspace{13.5pt}
Consider at first discrete time dynamical systems.
Let $M$ be a compact topological space with probability measure $\mu$ and
$f:M\to M$ be a homeomorphism. Consider some partition $\x$
(up to zero measure) $M=\sqcup\x_{i}$ of our space by
positive-measure sets with finite or countable number of indices.
We form the new partition $\xi^f_{-n}=\bigvee_{i=0}^nf^{-i}\x$, where
$\x\vee\e$ is a partition formed by the sets $\xi_\a\cap\e_\b$.

We define (\cite{KH}) entropy of the partition by
$H(\x)=-\sum\m(\x_j)\ln\m(\x_j)$ and entropy of a preserving measure
map $f$ w.r.t.\ $\x$ by
$h_\m(f;\x)=\lim\limits_{n\to\infty}\dfrac{H(\xi^f_{-n})}n$.
Then we define {\em measure entropy\/} of $f$ by
 $$
h_\mu(f)=\sup_{\{\x|h_\m(f;\x)<\infty\}} h_\m(f;\x).
 $$

To define {\em topological entropy\/} we should change partition $\x$ by
an open cover $\mathfrak U$ of $M=\cup U_j$, function $H(\x)$ by the
function $N(\mathfrak U)$ that is cardinality of minimal subcovering and
then take $\sup$ over all open coverings. This quantity is subject to the
following variational principle:
 $$
\hht(f)=\sup_\m h_\m(f),
 $$
where the supremum is taken over all $f$-invariant Borel probability
measures on $M$.

Another useful definition of the topological entropy goes for metric
spaces $(M,d)$. We define $d_n^f=\max\limits_{0\le i\le n-1}(f^i)^*d$ and
denote by $S(f,\ve,n)$ the minimal number of $\ve$-balls in
$d_n^f$-metric to cover $M$. Then
 \begin{equation}\label{topent}
\hht(f)
=\lim_{\ve\to0}\op{\overline{\lim}}\limits_{n\to\infty}
\dfrac{\ln S(f,\ve,n)}n
=\lim_{\ve\to0}\op{\underline{\lim}}\limits_{n\to\infty}
\dfrac{\ln S(f,\ve,n)}n.
 \end{equation}
This limit depends not on $d$ but on the $d$-topology only (\cite{KH}).

 \begin{examp}\po\label{ex90101}
Let $A\in\op{GL}_n(\Z)$. Then we can define map $A:T^n\to T^n$, where
$T^n=\R^n/\Z^n$. Let the measure be $d\mu=dx_1\we\dots\we dx_n$ in the
standard coordinates. Then
 \begin{equation}\label{90101}
\hht(A)=h_\mu(A)=\sum_{\l_A\in\op{Sp}^+(A)}\ln|\l_A|,
 \end{equation}
where $\op{Sp}^+(A)$ is the part of the spectrum outside the unit disk
$D_1\subset\CC$.
 \end{examp}

If $v\in\D(M)$ is a vector field on a manifold $M$ we define
$\hht(v)=\hht(\vp^1_v)$, where $\vp^t_v:M\to M$ is the flow of $v$. We
define the measure entropy similarly and note that for continuous time
dynamical systems the entropies are also connected by the variational
principle.

\subsection{\sl Compactness restriction.}

\hspace{13.5pt}
If $M$ is noncompact, the previous definitions basically do not work.
Still it is possible to define the entropy.

The easiest case is as follows. Suppose there exists a compact exhaustion
$M=\cup_{j=1}^\infty K_j$ with $K_j\subset K_{j+1}$ and each $K_j$ being
compact and $f$-invariant. Then we define topological entropy by
 $$
\hht(f)=\lim_{j\to\infty}\hht(\left.f\right|_{K_j})\in
\R_{\ge0}\cup\{+\infty\}.
 $$
Since the sequence $\hht(\left.f\right|_{K_j})$ is non-decreasing the limit
exists and since it exists for every exhaustion $\{K_j\}_1^\infty$ it does
not depend on its choice. We call such case
{\em non-compactness of the first type\/}.

As an example we consider Liouville integrable Hamiltonian system
$\op{sgrad}H$ with commuting integrals $I_1=H,I_2,\dots,I_n$ on $M^{2n}$,
for which generic intersection $\{I_j=\op{const}_j\}$ is a union of
tori. Then we obviously have invariant compact exhaustion and define
$\hht(\op{sgrad}H)$.

In the general situation we consider arbitrary compact exhaustion
$M=\cup K_j$ (not necessarily $f$-invariant). Let us study
the restriction of the metric $d_n^f$ to $K_j$.
Since $K_j$ is compact there is a finite covering by $\ve$-balls.
Calculating the asymptotic of the minimal cardinality of such a cover
we get similarly to (\ref{topent}) the quantity $\hht(f;K_j)$.
This sequence is non-decreasing by $j$ and we set
 $$
\hht(f)=\lim_{j\to\infty}\hht(f;K_j).
 $$
Again the limit does not depend on a choice of the compact exhaustion.
But it can be infinite even for sufficiently smooth $f$.
Let's call the situation {\em non-compactness of the second type\/}.

 \begin{examp}\po
Consider a linear automorphism $A:\R^n\to\R^n$. Then the entropy
$\hht(A)$ is given by the same formula (\ref{90101}) as in
example~\ref{ex90101}.
 \end{examp}

Note that this example shows $\hht(f^{-1})\ne\hht(f)$
generally in non-compact case. However all properties of the entropy
holds for non-compactness of the first type.

The measure analogs of the above definitions are straightforward.

\subsection{\sl Lyapunov exponents.}\label{S100101}

\hspace{13.5pt}
Let $\|\cdot\|$ be any norm on $TM$ and $f:M\to M$ be of class $C^1$.
Then we can consider the map
$v\mapsto\op{\overline{\lim}}\limits_{k\to\infty}\dfrac{\ln\|f^k_*v\|}k$.
This map takes values $\chi_1(x)\le\dots\le\chi_n(x)$ at (almost) any point
$x\in M^n$, which are called Lyapunov exponents (multiple values repeat).
They are defined similarly for the flows and do not depend on a choice of
the norm defining the topology.

Ruelle inequality says that if $\mu$ is $f$-invariant Borel probability
measure, then
 $$
h_\mu(f)\le\int_M\Biggl[\sum_{\chi_j(x)>0}\chi_j(x)\Biggr]d\mu.
 $$
In particular if the norm can be chosen so that $\|(\vp_t)_*v\|$ is
constant (Lebesgue) a.e., then the entropy vanishes w.r.t.\ any Liouville
(i.e.\ absolutely continuous) invariant probability measure.
Therefore we get

 \begin{examp}\po\label{ex:100101}
If a Hamiltonian system with Hamiltonian $H$ is completely integrable, then
for every Liouville invariant measure $h_\m(\op{sgrad}H)=0$.
 \end{examp}

Moreover Pesin theorem states that the inequality above becomes equality
 $$
h_\mu(f)=\int_M\sum_{\chi_j>0}\chi_j(x)d\mu
 $$
for $C^{1+\ve}$-diffeomorphisms/flows and any Liouville measure.

Since for integrable Hamiltonian system every measure
$\rho(I_1,\dots,I_n)|\O^n|$ is invariant, where $\O$ is the symplectic
form and $I_j$ are integrals, there are finite Liouville measures in this
case (even if $M^{2n}$ is non-compact).

\subsection{\sl Geodesic flows.}

\hspace{13.5pt}
Let $(M,g)$ be a compact Riemannian manifold. Let $\vp_t:TM\to TM$ be the
geodesic flow. Note that restriction of this flow to invariant compact
submanifold $T_cM=\{v\in TM,\|v\|_g=c\}$ is conjugated to
$\left.\vp_t\right|_{T_1M}$. Hence we define entropy $h(g)$
of the metric $g$ to be the entropy of the last flow.

Note that since the geodesic flow is Hamiltonian the (Liouville)
measure entropy vanishes $h_\m(g)=0$ for an integrable metric $g$ due to
example~\ref{ex:100101}. The behavior of topological entropy as we
will see in \S\ref{S4} can be more complicated.

 \begin{examp}\po
Geodesic flow on the group $\op{SO}(3)$ with left-invariant
Riemannian metric has $\hht(g)=0$.
 \end{examp}

Remark that canonical identification
$T^*M\stackrel{\sharp^g}\simeq TM$ allows
to consider the lifted flow $\vp_t:T^*M\to T^*M$. Its restriction to
the hypersurface $T^1M=\{H=1/2\}$ determines the entropy, where $H$ is the
Hamiltonian of the geodesic flow (for non-compact $M$ we should use the
above modification of $\hht$).

Consider now sub-Riemannian structure $(M,\Pi,g)$. We again can consider
geodesic flow and restrict it to the isoenergetic surface
$Q=\{H=1/2\}\subset T^*M$, $\vp_t:Q\to Q$. This $Q$ however is always
non-compact, so that we should use the non-compact version of the entropy:
 $\hht(\Pi,g):=\hht(\left.\vp_t\right|_{Q})$.

 \begin{examp}\po
For the group $\op{SO}(3)$ the topological entropy of
sub-Riemannian geodesic flow (see \S\ref{S3}) is zero.
This follows from the study of bifurcational diagram
(also directly from the description of the geodesics in \cite{VG}).
 \end{examp}

Note that if the Hamiltonian flow $\op{sgrad}H$ is Liouville-integrable
we have non-compactness of the first type, Lyapunov exponents vanish a.e.\
due to Liouville theorem and the (Liouville) measure entropy is zero.

\section{Commuting dynamical systems}\label{commut}

\hspace{13.5pt}
Let $f,g:M\to M$ be two commuting dynamical systems.
They define an action of $\Z^2$ on $M$. If we have more commuting systems
or commuting flows, then higher dimensional groups $\Z^k$ or $\R^k$ act.
There is a definition of entropy of such an action (\cite{S}). However if
it is non-trivial, then entropies of all dynamical systems given by
1-dimensional subgroups are infinite. We consider instead the case
when the higher entropy is zero and the entropies of the generators
of the group are finite.

 \begin{theorem}{\rm(\cite{Hu})\bf.}
Consider a $\Z^2$-action with generators $f,g$ of class $C^{1+\ve}$.
Let $h$ be either $h_\mu$, where $\mu$ is any $\Z^2$-invariant Borel
probability measure on $M$, or $\hht$ and then we additionally suppose
$fg\in C^\infty$. Then:
 \begin{equation}\label{eq:100101}
h(fg)\le h(f)+h(g).
 \end{equation}
 \end{theorem}

 \begin{proof}
Let us indicate the proof for a Liouville measure.
Let $\chi_i^f(x)$ and $\chi_j^g(x)$ be Lyapunov exponents for $f$ and $g$
respectively. Take $\chi_i^f(x)<\lambda<\chi_{i+1}^f(x)$. Then the map
$F_\l:=\l^{-1}f_*:TM\to TM$ commutes with $g_*$ and so
$(F_\l)^r\circ g_*=g_*\circ(F_\l)^r$. Tending $r\to\pm\infty$ we get
$g$-invariance of the repelling-expanding directions decomposition
$T_xM=H_i^-\oplus H_i^+$.

Applying this to all $\l$ different from Lyapunov exponents and
interchanging $f$ and $g$ we deduce for $x$ from a full measure
subset in $M$ the common decomposition $T_xM=\oplus H_s(x)$
which is: measurable, invariant and such that
on every subspace $H_s$ the Lyapunov exponents of both
$f$ and $g$ are constant. Therefore the Lyapunov exponents of $f\circ g$
are sums of the corresponding Lyapunov exponents for $f$ and $g$.

The claim now follows from the Pesin formula \ref{S2}\!\!\ref{S100101}.
The equality in (\ref{eq:100101}) is achieved if and only if we
sum positive exponents for $f$ and $g$ in the above decomposition
$\oplus H_i$ almost everywhere w.r.t.\ $\mu$.
 \end{proof}

The  case of general Borel measure as well as the inequality for $\hht$
is done in \cite{Hu} using Lyapunov charts (theorems B,C).

 \begin{rk}\po
Without commutativity assumption the statement is false. Actually the
Anosov automorphism $A:T^2\to T^2$  with
$A=\left(\begin{matrix} 2&1\\1&1\end{matrix}\right)$
can be decomposed $A=BB^t$ with
$B=\left(\begin{matrix} 1&1\\0&1\end{matrix}\right)$.
But $\hht(A)=\ln\dfrac{3+\sqrt5}2$, while $\hht(B)=\hht(B^t)=0$.

Also the requirement that the dynamical systems are smooth is crucial:
There are continuous $\Z^2$-actions (\cite{Pa}) with vanishing
directional entropies for all irrational directions and non-zero
for all rational.
 \end{rk}

 \begin{cor}\po\label{cor1}
Let $\Z^2$-dynamical system given by smooth commuting $f,g$
be either defined on a compact manifold or have non-compactness of the
first type. Let $h$ denote $h_\mu$ or $\hht$. If $h(g)=0$, then
$h(fg)=h(f)$.
 \end{cor}

 \begin{proof}
Actually $h(fg)\le h(f)$ and using $h(g^{-1})=h(g)=0$ we get the
reverse inequality.
 \end{proof}

Similar to (\ref{eq:100101}) formula holds for commuting flows
$\vp_\x^t,\vp_\e^s$ generated by vector fields $\x,\e$:
 \begin{equation}\label{eq:flows}
h(\vp^t_\x\vp^s_\e)\le h(\vp_\x^t)+h(\vp_\e^s),
 \end{equation}
where $h=\hht$ or $h_\mu$ as in the theorem.

Therefore if $\R^k$ acts on $M$ with generators
$\vp_1^{t_1},\dots,\vp_k^{t_k}$ we define
 \begin{equation}\label{pseudonorm}
\rho(v)=h(\vp_1^{t_1}\cdots\vp_k^{t_k})
\ \text{ for }\ v=(t_1,\dots,t_k)\in\R^k.
 \end{equation}

The inequality (\ref{eq:flows}) can be seen now as triangle inequality
and hence $\rho$ is a pseudonorm on $\R^k$ (provided non-compactness is
no greater than the first type). This $\rho$ however should not be
continuous and can be degenerate (not norm).

\section{Integrability of sub-Riemannian metrics on 3-dimensional Lie
groups}\label{S3}

\hspace{13.5pt}
Consider 3-dimensional Lie groups $G$. Every left-invariant structure
is presented on the Lie algebra $\G$. So we write the classification of
left-invariant non-holonomic metrics on contact structures along Bianchi
classification of 3-dimensional Lie algebras (\cite{VG}). We identify
proportional metrics (the geodesic flows are reparame\-te\-rized). We
present distribution $\Pi\subset\G$ and the metric by an orthonormal basis
$\x_1,\x_2$:
 \begin{enumerate}
  \item Heisenberg algebra ${\mathfrak h}(3)$: $[e_1,e_2]=e_3,[e_k,e_3]=0$
   and $\x_k=e_k$, $k=1,2$.
  \item Solvable algebras: $[e_1,e_2]=0$,
   \begin{enumerate}
  \item $[e_1,e_3]=\l_1e_1,[e_2,e_3]=\l_2e_2$ ($\l_1\ne\l_2$)
   and $\x_1=e_1+e_2,\x_2=e_3$.
  \item $[(e_1+ie_2),e_3)]=e^{-i\vp}(e_1+ie_2)$ ($0<\vp<\pi$)
   and $\x_1=e_1,\x_2=e_3$.
  \item $[e_1,e_3]=e_1+e_2,[e_2,e_3]=e_2$
   and $\x_1=e_1,\x_2=e_3$.
   \end{enumerate}
  \item Orthogonal algebra $\op{so}(3)$:
 $[e_1,e_2]=e_3,[e_2,e_3]=e_1,[e_3,e_1]=e_2$ and $\x_1=e_1,\x_2=\z e_2$.
  \item Special linear algebra $\op{sl}_2(\R)$:
  $[e_1,e_2]=e_3,[e_1,e_3]=2e_1,[e_2,e_3]=-2e_2$ and two possibilities
    \begin{enumerate}
  \item $\x_1=e_1,\x_2=\z e_2 $;
  \item $\x_1=e_1+e_2,\x_2=\z e_3$.
   \end{enumerate}
  In three last cases the parameter $\z\in\R_+$.
 \end{enumerate}
One can integrate the equations of geodesics (see \cite{VG} for a
description in terms of semi-direct product). We consider Liouville
integrability.

 \begin{theorem}\po
Non-holonomic geodesic flows of the above metrics on 3-di\-men\-si\-o\-nal
Lie groups are $C^\oo$-integrable.
 \end{theorem}

 \begin{proof}
First note that any vector $v\in\G$ generates left- and right-invariant
vector fields $L_v$ and $R_v$ on $G$, which can be considered as functions
on $T^*G$. We also denote by $L_v,R_v\in C^\infty(T^*G)$ the invariant
extensions of any function $v\in C^\infty(\G^*)$, not necessarily linear.
Any two $L_v$ and $R_w$ Poisson commute and the other Poisson brackets are
$\{L_v,L_w\}=L_{[v,w]}$, $\{R_v,R_w\}=R_{[v,w]}$ (\cite{AKN}, \cite{F}),
where the bracket $[,]$ on $C^\infty(\G^*)$ is induced by the usual
commutator (this bracket is called Lie-Poisson or
Berezin-Kirillov-Konstant-Souriau).

So finding one left-invariant function $I_2=L_v$ commuting with
$2H=\x_1^2+\x_2^2$ will suffice, just take any $I_3=R_w$.
We take the Casimir function $F$ for $v$, i.e.\ such a function that
$[\cdot,F]\equiv0$. This function always exists locally since $\G$ has odd
dimension, but globally it can have singularities.

In the case of $\op{so}(3)$ we have smooth Casimir function
$F=e_1^2+e_2^2+e_3^2$. For $\op{sl}_2(\R)$ $F=4e_1e_2-e_3^2$. For the
Heisenberg group $F=e_3$. For the solvable group we consider $[\cdot,e_3]$
as a linear vector field on the plane $\R^2(e_1,e_2)$. We take $F$ to be
an integral of this vector field. For example for the case (3a) we take
$F=(e_1)^{\l_2}(e_2)^{-\l_1}$.

Note that these functions $F$ have singularities along some axes in the
solvable case. To overcome this difficulty we use the lemma.
 \begin{lem}\po
If 3-dim Lie group $G$ is not semisimple, then for every left-invariant
$H$ there exist two commuting right-invariant integrals $I_2,I_3$.
 \end{lem}

Actually due to Bianchi classification we can always find a
two-dimensional commutative subalgebra $\langle v,w\rangle\subset\G$ in
non-semisimple case ($\langle e_1,e_2\rangle$ for solvable algebras). Now
it's easy to check that in every case of non-holonomic metric we obtain 3
functionally independent a.e.\ integrals.
 \end{proof}

{\bf Note on the non-holonomic flows on $\op{SL}_2(\R)$.}
The two contact structures on $\op{SL}_2(\Z)$ have nice geometric
interpretations (\cite{Kr}). It is well known that
$\op{SL}_2(\Z)/{\{\pm1\}}\simeq ST^*L^2$ for the spherical bundle of the
Lobachevskii plane $L^2$ and since every surface $M^2$
of genus $g>1$ can be obtained as quotient of $L^2$ by a discrete
subgroup we end up with two contact structures on $ST^*M^2$ of
which one is the standard and the other is the connection
form associated with a metric of constant negative curvature.

Note that the metric descends to $ST^*M^2$ and for the second structure
even to $M^2$. However we lose integrals: On $ST^*L^2$ we have all two
additional integrals, on $ST^*M^2$ we loose one right-invariant and under
projection to $M^2$ we loose the other integral. This is consistent with
the Kozlov theorem \cite{Ko} according to which the geodesic flow on
$M^2$ with $g>1$ handles is not analytically integrable. Still we see
there is an integrable lift of the flow.

\section{Non-holonomic metrics on suspensions $M_A$}\label{S4}

\subsection{\it Construction of sub-Riemannian structure.}

\hspace{13.5pt}
Let $M_A=T^2\tilde\t_AS^1$ be quotient of the
cylinder $\C=T^2(\vp_1,\vp_2)\t\R^1(\vp_3)$ w.r.t.\ the free action given
by the map $\hat A(\vp_1,\vp_2,\vp_3)=(A(\vp_1,\vp_2),\vp_3+1)$.

\bff 1
Consider at first the case of semisimple $A\in\op{SL}_2(\Z)$
with real eigenvalues $\l>1,\l^{-1}\in(0,1)$. Let
$\e_1,\e_2\in T_*(T^2)$ be the eigenvectors
(the case $\l=1$ was considered in example~\ref{ex:xe} and
the case $\l<0$ will follow similarly because the distribution $\Pi$
has no canonical orientation and so the change
$\e_j\mapsto-\e_j$ preserves the metrics below).
The action on the basis is the following
$\hat A:(\e_1,\e_2,\p_{\vp_3})\mapsto(\l\e_1,\l^{-1}\e_2,\p_{\vp_3})$.

We define the contact structure on $\C$ by the vector fields
$\x_1=e^{-\ln\l\cdot\vp_3}\e_1+e^{\ln\l\cdot\vp_3}\e_2$, $\x_2=\p_{\vp_3}$:
$\Pi=\langle\x_1,\x_2\rangle$. Moreover we fix non-holonomic metric $g$ by
requiring that $\x_1,\x_2$ is an orthonormal basis.

 \begin{lem}\po
$\Pi$ and $g$ are invariant w.r.t\ the action $\hat A$.\qed
 \end{lem}

 \begin{lem}\po
The distribution $\Pi$ is contact.
 \end{lem}

 \begin{proof}
 $$
[\x_2,\x_1]=\x_3=
-\ln\l e^{-\ln\l\cdot\vp_3}\e_1+\ln\l e^{\ln\l\cdot\vp_3}\e_2.
 $$
But $\x_3\notin\Pi$ because the matrix
 $$
\left[\begin{matrix}\x_1\\ \x_3\end{matrix}\right]=
\left(\begin{matrix}
e^{-\ln\l\cdot\vp_3} & e^{\ln\l\cdot\vp_3}\\
-\ln\l e^{-\ln\l\cdot\vp_3} & \ln\l e^{\ln\l\cdot\vp_3}
\end{matrix}\right)
 $$
is Vandermonde and hence nondegenerate for $\ln\l\ne0$.
 \end{proof}

Summarizing we get a sub-Riemannian structure $(\Pi,g)$ on $M_A$.

\bff 2
Now let the eigenvalues of $A$ be $\l,\bar\l=e^{\pm\te_k}$. Since
the matrix is integer-valued we have $\te_k=\pm\dfrac{2\pi}k$, $k=3,4,6$
(modulo the cases we have considered). Thus $A$ has the matrix
$\left(\begin{matrix} \cos\te_k & \sin\te_k \\
-\sin\te_k & \cos\te_k \end{matrix}\right)$ in a basis
$\e_1,\e_2$ of $T_*(T^2)$. We define
$\x_1=\cos(\te_k\vp_3)\e_1-\sin(\te_k\vp_3)\e_2$, $\x_2=\p_{\vp_3}$
and as before $\Pi=\langle\x_1,\x_2\rangle$, $g=(\xi_1^*)^2+(\x_2^*)^2$.
One easily checks that $(\Pi,g)$ is invariant, non-holonomic and hence
defines a sub-Riemannian structure on $M_A$.

\bff 3
Consider finally the Jordan box, i.e.\ $A$ is conjugated to
$\left(\begin{matrix} 1&1\\0&1\end{matrix}\right)$.
Let $\e_1,\e_2$ be the corresponding basis,
$\hat A:(\e_1,\e_2)\mapsto(\e_1,\e_1+\e_2)$. We set
$\x_1=\cos(2\pi\vp_3)\e_1+\sin(2\pi\vp_3)\bigl(\e_2-\vp_3\e_1\bigr)$,
$\x_2=\p_{\vp_3}$ and proceed as before to get sub-Riemannian structure
$(\Pi,g)$ on $M_A$.

\subsection{\it $C^\infty$-integrability.}

\hspace{13.5pt}
The Hamiltonian of the sub-Riemannian flow is $H=\frac12[\x_1^2+\x_2^2]$.
Let $p',p'',p_3$ be momenta dual to vectors $\e_1,\e_2,\p_{\vp_3}$ from
$\C$ (i.e.\ evaluation matrix $p(\eta)$ is $E_{3\t3}$).

This Hamiltonian is $\hat A$-invariant, so that it descends to
$M_A$. Since $H$ does not depend on $\vp_1,\vp_2$ the functions
$p_1,p_2$ and hence $p',p''$ are integrals on $\C$.
However these functions are not $\hat A$-invariant.
We get invariants from their combination.

\bff 1
In this case we can write
 \begin{equation}\label{Hamilt}
2H=(e^{\ln\l\cdot\vp_3}p'+e^{-\ln\l\cdot\vp_3}p'')^2+(p_3)^2.
 \end{equation}
Here the functions transform by the rule
${\hat A}:(p',p'',p_3)\mapsto(\l^{-1}p',\l p'',p_3)$.
Thus we have the following integrals: $I_1=H$ and
 \begin{equation}\label{BT-in1}
I_2=p'p'',\quad
I_3=\sin\Bigl(2\pi\frac{\ln|p'|}{\ln\l}\Bigr)e^{-(p'p'')^{-2}},
 \end{equation}
which are invariant and functionally independent a.e.
Thus we conclude $C^\infty$-integrability of
sub-Riemannian structure $(\Pi,g)$ on $M_A$.

\bff 2
In this case the action is the rotation by $\te_k$ and so the additional
integrals are:
 \begin{equation}
I_2=(p')^2+(p'')^2,\quad
I_3=\op{Re}(p'+ip'')^k.
 \end{equation}
Note that in this case the integrals are analytic.

\bff 3
The action is ${\hat A}:(p',p'',p_3)\mapsto(p'-p'',p'',p_3)$. Therefore
the additional integrals are:
 $$
I_2=p'',\quad
I_3=\sin\Bigl(2\pi\frac{p'}{p''}\Bigr)e^{-(p'')^{-2}}.
 $$
This is again the case of $C^\infty$, not $C^\oo$-integrals.

\subsection{\it Obstructions to  "nice" integrability.}

\hspace{13.5pt}
Here we generalize Taimanov's result.

 \begin{dfn}
We say sub-Riemannian structure $(M^m,\Pi,g)$
is {\em geometrically simple\/} if $Q=\{H=1/2\}$ contains a closed nowhere
dense invariant subset $\Gamma$ such that:
 \begin{enumerate}
  \item $Q\setminus\Gamma=\bigcup_{j=1}^d U_j$, where $U_j$ are open and
path-connected.
  \item Each $U_j$ is foliated by Liouville tori,
$U_j\simeq T^{m}\t D^{m-1}$.
  \item[\it 3\,$'$\!\!.]
Let $\overline Q$ be the one-point compactification along fibers
of the projection $p:Q\to M$. Then for every $x\in\overline Q$
there exists an arbitrary small neighborhood $W_x$ such that
$W_x\setminus\Gamma$ has a finite number of path-connected components.
 \end{enumerate}
 \end{dfn}

The condition $3'$) is a version of the corresponding 3) from~\cite{T}
for non-compact $Q$. In fact both conditions are equivalent to the
following:
 \begin{enumerate}
  \item[\it 3\,$''$\!\!.]{\it
For every $x\in M$ there exist arbitrary small
neighborhoods $U_x\subset V_x$ such that
$|\pi_0\bigl(\bigl[p^{-1}(V_x\setminus{\overline U}_x)\bigr]\cup
\bigl[p^{-1}({\overline U}_x)\setminus\Gamma\bigr]\bigr)|<\infty$.
}
 \end{enumerate}
We remark also that the set $\Gamma$ is usually bigger than the
bifurcational set $\Sigma\subset Q$ in Liouville integrable case.

 \begin{theorem}\po
Let a sub-Riemannian structure $(\Pi,g)$ on a closed manifold $M^m$ be
Liouville integrable and geometrically simple. Then:
 \begin{itemize}
  \item The group $\pi_1(M)$ is almost commutative, i.e.\ it contains a
commutative subgroup of finite index.
  \item The first Betti number $\dim H_1(M^m;\R)\le m$.
 \end{itemize}
 \end{theorem}

 \begin{proof}
Actually one take a graph with vertices numerated by the components of
$p^{-1}(V_x\setminus{\overline U}_x)\cup
\bigl(p^{-1}({\overline U}_x)\setminus\Gamma\bigr)$.
This number is finite due to $3''$).
Applying sub-Riemannian Hopf-Rinow theorem \cite{BR} we
get labeled directed edges of this graph with labels $s\in\pi_1(M,x)$
realized as geodesic loops lifted to
$p^{-1}(V_x)\cap(Q\setminus\Gamma)$.
The rest of Taimanov's proof \cite{T} goes now unchanged.
 \end{proof}
\smallskip

We call the set of a.e.\ independent involutive integrals
$I_1=H,I_2,\dots,I_n$ of Hamiltonian system on $W^{2n}$
{\em geometrically simple\/} on $Q_c=\{H=c\}$ if the
bifurcation sets $\Sigma$ are diffeomorphic on $Q_{c'}$ for close
$c'\approx c$ and the stratification of $Q_c$ by singularities
of the moment map $Q_c\ni x\mapsto(I_2(x)\dots,I_n(x))$ can
be extended to stratification $Q_c=\cup{\tilde\Sigma}^\a_j$ with
$\Gamma=\cup_{j>0}{\tilde\Sigma}^\a_j$ giving the conditions 1-3) of the
previous definition. Certainly if sub-Riemannian flow possesses a
geometrically simple involutive set of integrals, it is geometrically
simple.

For compact $Q$ (Riemannian case) analytic involutive set of integrals
is geometrically simple \cite{T}. In sub-Riemannian case
non-compactness makes difficulties for application of the standard
theorems about sub-analytic (constructive) sets. Yet known integrable
cases are geometrically simple. Actually if the integrals are $C^\oo$ and
polynomial by momenta (in this case we have Liouville tori),
the set $(I_j)_1^n$ is geometrically simple.
This is because singularities of $(I_j)$ give stratification of $Q$ at
infinity, while the finite (compact) part can be additionally
substratified (and even made into simplical decomposition as in \cite{T}).

 \begin{cor}\po
Since $\pi_1(M_A)$ is not almost commutative for $A$ with real eigenvalues
different from $\pm1$ and for the Jordan box~\cite{GK}, the
sub-Riemannian geodesic flow is not Liouville integrable with
geometrically simple set of integrals. \qed
 \end{cor}

\subsection{\it Calculation of topological entropy.}

\hspace{13.5pt}
Calculations of topological entropy in the case of rotation and Jordan
box are easy. So we consider only semi-simple case.

First note that the following submanifold $\{p'=0,p''=0,p_3=1\}$ of
$T^*\C$ is invariant under the transformation $\hat A$. So this
3-dimensional submanifold descends to $T^*M_A$ and belongs to the
hypersurface $\{H=1/2\}$. Moreover it is invariant under the geodesic flow
and the flow moves with the unit speed along $\p_{\vp_3}$. So fixing
$\vp_3=0$ we get the torus $T^2$ and the Poincar\'e return map of the
geodesic flow induces the "Anosov map" $T^2\to T^2$ with matrix $A$.
Since the entropy of a system is not less than entropy of any subsystem,
we conclude: $\hht(\op{sgrad}H)\ge\ln|\op{Sp}^+(A)|$.

Now let us prove the inverse. For this we study bifurcations of
$(I_1,I_2,I_3)$. Let us note that the integrals (\ref{BT-in1}) coincide
with the integrals of Bolsinov-Taimanov \cite{BT2}, while the Hamiltonians
are different. The bifurcation set in the Riemannian case is twice bigger
than in the sub-Riemannian one, but still the systems behave similarly.

It is easy to see that $\op{rk}\{dI_1,dI_2,dI_3\}=3$ along
a trajectory if $(p')^2+(p'')^2\ne0$ and either $p_3\ne0$ or
$e^{\ln\l\cdot\vp_3}p'\ne e^{-\ln\l\cdot\vp_3}p''$ at least at one
(and then any) point of it. Then by the Liouville theorem this trajectory
lies on an invariant torus and the motion is quasi-periodic.

Let $p_3=0\Rightarrow \vp_3=\op{const}$ and $p''=e^{2\ln\l\cdot\vp_3}p'$,
but $(p')^2+(p'')^2\ne0$. Then the point moves quasi-periodically
along the fiber-torus $T^2\subset M_A$ and hence again contributes nothing
to the topological entropy.

Now on the invariant set $N=\{p'p''=0\}$ our Hamiltonian (\ref{Hamilt})
coincides with the one from \cite{BT2}. Therefore the description of the
trajectories is the same. In particular $N=N'\cup N''$, each of the
summands $N'=\{p'=0\}$, $N''=\{p''=0\}$ being invariant and
diffeomorphic to $M_A\t S^1$. The intersection
$N'\cap N''$ has two components $V',V''\simeq M_A$ given by
$\{p_3=\pm1\}$ such that $V'$ is a stable manifold for $N'$ and unstable
for $N''$, $V''$ is stable for $N''$ and unstable for $N'$.

Moreover the Hamilton equation of the geodesic flow implies that Lyapunov
exponents on $V$ are exactly $\pm\ln\l,0$. So by Katok theorem the only
invariant Borel measure on $N$ can be supported on $V'\cup V''$. But the
flow on $V'$ and $V''$ is Anosov with Poincar\'e time-one map
$A:T^2\to T^2$. Therefore the variational principle yields $\hht=\ln\l$.

\subsection{\it Reeb vector field.}

\hspace{13.5pt}
We prove in this section theorem~\ref{th2}.

The unique (up to sign) contact form $\a$ for sub-Riemannian structure is
given by the equality $d\a(\x_1,\x_2)=-1$ or $\a([\x_1,\x_2])=1$, where
$\x_1,\x_2$ is an arbitrary $g$-orthogonal frame from $C^\infty(\Pi)$.

\bff 1 Consider at first the semisimple case
$\x_1=e^{-\ln\l\cdot\vp_3}\e_1+e^{\ln\l\cdot\vp_3}\e_2$, $\x_2=\p_{\vp_3}$.
The Reeb field satisfies $\nu_\a=v\op{mod}\Pi$ for
$v=[\x_1,\x_2]=
\ln\l\bigl(e^{-\ln\l\cdot\vp_3}\e_1-e^{\ln\l\cdot\vp_3}\e_2\bigr)$.

 \begin{lem}\po
$v$ is a symmetry of $\Pi$.
 \end{lem}

 \begin{proof}
Actually $[v,\x_1]=0$, $[v,\x_2]=\ln^2\l\x_1$.
 \end{proof}

 \begin{cor}\po
$\nu_\a=v$ is the Reeb vector field.
 \end{cor}

 \begin{proof}
Actually $\a(v)=1$.
 \end{proof}

Now we see that the flow of $\nu_\a$ is quasiperiodic on the tori-fibers
of $M_A$. Therefore $\hht(\nu_\a)$ vanish.

\bff 2
The case of rotation is absolutely similar.

\bff 3
For the Jordan box the field $v=[\x_1,\x_2]$ is not a symmetry of $\Pi$
and hence is not the Reeb field. Calculation shows that the contact form
is $\a=\bigl(2\pi+\sin^2(2\pi\vp_3)\bigr)^{-1}
\bigl(\sin(2\pi\vp_3)\e_1^*+(\vp_3\sin(2\pi\vp_3)-\cos(2\pi\vp_3)\e_2^*)
\bigr)$,
where $\e_1^*,\e_2^*$ is the co-basis of $T^*T^2$ dual to $\e_1,\e_2$.
Now since $d\a=d\vp_3\we\p_{\vp_3}(\a)$, the Reeb field is
$\nu_\a=v+b(\vp_3)\x_1$ and again $\hht(\nu_\a)=0$.

 \begin{rk}\po
Though $\hht(\hat g)=\hht(g)$ in theorem~\ref{th3}, it is not a simple
combination of corollary~\ref{cor1} and theorem~\ref{th1}:
formula (\ref{eq:100101}) is not working here.
 \end{rk}

\section{Reeb-symmetric SR-structures}

\hspace{13.5pt}
Let $\Pi$ be a contact structure on $M^{2n+1}$ and $\nu_\a$ be the Reeb
vector field as in introduction.

 \begin{dfn}
Sub-Riemannian structure $(\Pi,g)$ is called {\em Reeb-symmetric\/} if
$\nu_\a$ is a symmetry of $g$ on $\Pi$: $L_{\nu_\a}g=0$. In other words
$\nu_\a$ is a Killing vector field.
 \end{dfn}

 For instance the Heisenberg SR-structure (example~\ref{HSR}) is of this
kind.
\newline
 Let us define the function $I_g(p)=\langle p,\nu_\a\rangle^2$ (or just
$\nu_\a^2$) on $T^*M$. We will call it {\em Reeb momentum\/}.

 \begin{prop}\po
$(\Pi,g)$ is Reeb symmetric iff the momentum $I_g$ is the first integral
of the sub-Riemannian geodesic flow.
 \end{prop}

 \begin{proof}
Actually locally in Darboux coordinates: $\a=dx-\sum_{i=1}^n z_idy_i$ and
$\nu_\a=\p_x$.  Now it's obvious that $L_{\p_x}g=0$ and $\{I_g,H\}=0$ are
equivalent to the condition that $g$ does not depend on $x$
[lifted from some $g_0$ on $\R^{2n}(y,z)$].
 \end{proof}

Let us construct by $g$ and $\nu_\a$ the Riemannian metric $\hat g$ on $M$
as in introduction. Then the Riemannian and sub-Riemannian Hamiltonians
are connected by the formula
 \begin{equation}\label{H+H}
H_{\hat g}=H_g+I_g.
 \end{equation}

 \begin{cor}\po
Reeb momentum $I_g$ is an integral of Riemannian geodesic flow if and only
if it is an integral for sub-Riemannian one. \qed
 \end{cor}

Therefore one can integrate instead of $g$, the metric $\hat g$, finding
involutive set of integrals starting from $I_1=H_{\hat g}$ and $I_2=I_g$.

 \begin{rk}\po
Let us note that the metrics ${\hat g}^{(t)}$, which equal $g$ on
$\Pi$ and $1/{\sqrt t}$ on orthogonal $\nu_\a$, satisfy:
$\lim\limits_{t\to +0}\op{dist}_{{\hat g}^{(t)}}=\op{dist}_g$ --
sub-Riemannian distance. The corresponding Hamiltonians
$H_{{\hat g}^{(t)}}=H_g+tI_g$.
 \end{rk}

Let us note also about topological entropy:

 \begin{prop}\po
If $(\Pi,g)$ is smooth and Reeb-symmetric, then vanishing of two of three
$\hht(\op{sgrad}H_{\hat g})$, $\hht(\op{sgrad}H_g)$, $\hht(\nu_\a)$
implies vanishing of the third. Moreover if $\hht(\nu_\a)=0$, then
$\hht(g)=\hht(\hat g)$.
 \end{prop}

 \begin{proof}
Actually this follows from (\ref{H+H}) and commutativity.
 \end{proof}

 \begin{rk}\po
Note that even though the Reeb field $\nu_\a$ is {\em divergence-free\/}
and shares some good properties (no attractors/repellers), we can have
$\hht(\nu_\a)>0$. Examples are provided by geodesic flows on negatively
curved manifolds.
 \end{rk}

\section{Concluding remark: norm via entropy}\label{S6}

\hspace{13.5pt}
Any integrable Hamiltonian system on $(M^{2n},\oo)$ with integrals
$(I_1,\dots,I_n)$ determines $\R^n$-action with generators $\vp_j^{t_j}$
being shifts along $\op{sgrad}_\oo(I_j)$. Then the formula
(\ref{pseudonorm}) with $h=\hht$ gives a pseudonorm
$\rho=\rho_{\text{top}}$ on this $\R^n$ (of course for $h=h_\mu$ with
Liouville invariant measure $\mu$ we have $\rho_\mu\equiv0$).

For instance if the action reduces to the action of the torus $T^n$ then
$\rho_{\text{top}}\equiv0$. This is obvious from the definition (also
is seen from Atiyah convexity theorem).

Examples of this and \cite{BT2} papers provide $\R^3$-action with
$\rho_{\text{top}}\not\equiv0$. However $\rho_{\text{top}}$ is not
a norm, because $\hht(\op{sgrad}I)=0$ for any $I=\l_2I_2+\l_3I_3$.

 \begin{theorem}\po
There is an integrable Hamiltonian system with
nondegenerate $\rho_{\text{top}}$ (norm).
 \end{theorem}

To explain the result consider a piece of bifurcation diagram for the
momentum map $(I_1,I_2,I_3)$ in the semisimple case that is pictured in
Fig.~1. The central vertical line corresponds to the most complicated
singularity $\{p_3=\pm1\}$ and in the preimage of each of its point we
have 2 copies of $M_A$. Two curved boundary surfaces correspond to
bifurcation diagram either, while 4 other plane pieces of the boundary
consist of regular points and appear due to the size of a neighborhood.

 \begin{picture}(500,147)
\put(-130,150){\centerline{{\special{em:graph bifdgr01.gif}}}}
\put(140,10){\bf Fig.\ 1}
 \end{picture}

Now let's take 3 mutually non-intersecting such pieces in the affine
$\R^3$ such that their central axes are pair-wise non-parallel (we can
direct them along axes $I_1,I_2,I_3$). We connect their regular points
by small tubes as shown in Fig.~2.

 \begin{picture}(500,149)
\put(-140,150){\centerline{{\special{em:graph bifdgr02.gif}}}}
\put(160,10){\bf Fig.\ 2}
 \end{picture}
This is equivalent to connecting
preimages of pieces by the tubes $T^3\t I\t D^2(\ve)$. The tubes in $\R^3$
carry the affine structure, that permits by Arnold-Liouville theorem to
define symplectic form and a Hamiltonian on the obtained manifold which is
integrable with the Hamiltonian action of $\R^3$ given by $(I_1,I_2,I_3)$
and such that $\rho_{\text{top}}$ is a norm on this $\R^3$.

 \pagebreak

\vspace{-6pt}
\hspace{-20pt} {\hbox to 12cm{ \hrulefill }}

{\footnotesize
\hspace{-10pt}
Inst. of Math., University of Tromsoe, Tromsoe 90-37, Norway;\quad
kruglikov\verb"@"math.uit.no}

\end{document}